\documentclass[%
  a4paper,
  twocolumn,
]{mypreprint}

\usepackage[english]{babel}

\usepackage{graphicx}

\usepackage{amsmath}
\usepackage{amssymb}
\usepackage{amsthm}

\usepackage{cleveref}

\providecommand\xSDRE[1]{\texttt{xSDRE-#1}}
\providecommand\LQR{\texttt{LQR}}
\providecommand\tc{t_{\mathsf{c}}}

\newcommand{\R}{\ensuremath{\mathbb{R}}}
\newcommand{\N}{\ensuremath{\mathbb{N}}}

\newcommand{\tA}{\ensuremath{\widetilde{A}}}
\newcommand{\tB}{\ensuremath{\widetilde{B}}}
\newcommand{\tN}{\ensuremath{\widetilde{N}}}
\newcommand{\tv}{\ensuremath{\tilde{v}}}
\newcommand{\hA}{\ensuremath{\widehat{A}}}
\newcommand{\hN}{\ensuremath{\widehat{N}}}
\newcommand{\hr}{\ensuremath{\hat{r}}}
\newcommand{\hv}{\ensuremath{\hat{v}}}
\newcommand{\hy}{\ensuremath{\hat{y}}}
\newcommand{\hrho}{\ensuremath{\hat{\rho}}}
\newcommand{\hmu}{\ensuremath{\hat{\mu}}}
\newcommand{\hnu}{\ensuremath{\hat{\nu}}}

\newcommand{\Azcl}{\ensuremath{A_{0,\mathsf{cl}}}}

\newcommand{\Nl}{\ensuremath{N_{\lambda}}}

\newcommand{\trans}{\ensuremath{\mkern-1.5mu\mathsf{T}}}

\begin{document}


\title{Low-complexity linear parameter-varying approximations of
  incompressible Navier-Stokes equations for truncated state-dependent Riccati
  feedback}
  
\author[$\ast$]{Jan Heiland}
\affil[$\ast$]{Max Planck Institute for Dynamics of Complex Technical
  Systems, Sandtorstra{\ss}e 1, 39106 Magdeburg, Germany.
  \email{heiland@mpi-magdeburg.mpg.de}, \orcid{0000-0003-0228-8522}
  \authorcr \itshape
  Otto von Guericke University Magdeburg, Faculty of Mathematics,
  Universit{\"a}tsplatz 2, 39106 Magdeburg, Germany.
  \email{jan.heiland@ovgu.de}}

\author[$\dagger$]{Steffen W. R. Werner}
\affil[$\dagger$]{Courant Institute of Mathematical Sciences, New York
  University, New York, NY 10012, USA.\authorcr
  \email{steffen.werner@nyu.edu}, \orcid{0000-0003-1667-4862}}
  
\shorttitle{Low-complexity LPV approximations of NSEs for SDRE feedback}
\shortauthor{J. Heiland, S. W. R. Werner}
\shortdate{2023-03-21}
\shortinstitute{}

\keywords{}

\msc{}

\abstract{%
  Nonlinear feedback design via state-dependent Riccati equations is well
  established but unfeasible for large-scale systems because of computational
  costs. 
  If the system can be embedded in the class of linear parameter-varying (LPV)
  systems with the parameter dependency being affine-linear, then the
  nonlinear feedback law has a series expansion with constant and precomputable
  coefficients.
  In this work, we propose a general method to approximating nonlinear systems
  such that the series expansion is possible and efficient even for
  high-dimensional systems.
  We lay out the stabilization of incompressible Navier-Stokes equations as
  application, discuss the numerical solution of the involved matrix-valued
  equations, and confirm the performance of the approach in a numerical example.
}

\novelty{}

\maketitle

  
\section{Introduction}%
\label{sec:intro}

Nonlinear feedback design for large-scale systems is challenging, as both the
complexity induced by nonlinearities and the huge computational tasks caused by
the system's size have to be resolved. 
The commonly used methods of \emph{backstepping}~\cite{Kok92}, feedback
linearization~\cite[Ch. 5.3]{Son98}, or \emph{sliding mode}
control~\cite{Dod15} require structural assumptions and, thus, may not be
accessible to a general computational framework. 
The both holistic and general approach via the Hamilton-Jacobi-Bellman (HJB)
equations, however, is only feasible for very moderate system sizes or calls
for model order reduction; see, e.g.,~\cite{BreKP19} for a relevant discussion
and an application in fluid flow control. 
As an alternative to reducing the system's size, one may consider approximations
to the solution of the HJB equations of lower complexity.
For that, for example, truncated polynomial expansions~\cite{BreKP18} or 
suboptimal solutions via the so called state-dependent Riccati equation
(SDRE)~\cite{BanLT07} are considered.
Here, we will follow on recent developments~\cite{AllKS23} on series expansions
of the SDRE approximation to the HJB solution that can mitigate the still high
computational costs of repeatedly solving high-dimensional Riccati equations.

As the general setup, we consider the input-affine system
\begin{equation} \label{eq:gen-affine-system}
  \dot v(t) = f(v(t)) + Bu(t), \quad y(t) = C v(t),
\end{equation}
where for time $t > 0$, $x(t)\in \R^{n}$ denotes the state, $u(t)
\in \R^{p}$ and $y(t)\in \R^{q}$ denote the input
and output, $f\colon \R^{n}\to \R^{n}$ is a possibly nonlinear function,
and where $B$ and $C$ are linear input and output operators.
Under the mild condition that $f$ is Lipshitz continuous and $f(0) = 0$, one can
factorize the nonlinearity  $f(v) = A(v) v$ with a state-dependent coefficient
matrix $A(v)$ and bring the system~\cref{eq:gen-affine-system} into 
state-dependent coefficient (SDC) form:
\begin{align}\label{eq:gen-sdc-sys}
  \dot v(t) = A(v(t)) v(t) + Bu(t), \quad y(t) = Cv(t);
\end{align}
see, e.g.,~\cite[Eq.~(7)]{BenH18}.
For such systems, one can define a feedback by
\begin{equation*}
  u(t) = -B^{\trans} P(v(t)) v(t),
\end{equation*}
where $P$ solves the SDRE
\begin{equation}\label{eq:sdre}
  A(v)^{\trans} P+PA(v) - PBB^{\trans} P = -C^{\trans} C;
\end{equation}
see~\cite{BanLT07} for general principles and~\cite{BenH18} for a proof of
performance beyond an asymptotic smallness condition.
Because of its nonlinear and, possibly high-dimensional nature, a
solve of the SDRE~\cref{eq:sdre} comes at high costs that make the
SDRE approach unfeasible for large systems; see~\cite{BenH18} for an
example illustrating how the effort grows with the system's dimension.

If, however, the factorization $f(v) = A(v) v$ is affine-linear with respect to
a parametrization $\rho(v) \in \R^{m}$ of $v$, i.e., it can be
represented as
\begin{equation}\label{eq:affine-LPV-representation}
  A(v) = A_{0} + \sum_{k=1}^{m} \rho_{k}(v) A_{k},
\end{equation}
then the solution $P(v)$ of the SDRE has a \emph{first-order} approximation
\begin{equation*}
  P(v) \approx P_{0} + \sum_{k=1}^{m} \rho_k(v)L_{k}
\end{equation*}
where $P_{0}$ and $L_{k}$, $k = 1, \ldots, m$ can be precomputed by one Riccati
and $m$ Lyapunov equations; see~\cite{AllKS23, BeeTB00}.
In this work we propose a general approach for controller design that bases on
approximative representations as in~\cref{eq:affine-LPV-representation} for
which we employ
\begin{itemize}
  \item an SDC representation of the nonlinear system,
  \item an approximative parametrization $\hrho(v)\in \R^{r}$, where
    approximative means that $r$ is chosen small such that an exact
    reconstruction of $v$ from $\hrho(v)$ might not be possible, and
  \item an affine-linear approximation of the
    coefficient~\cref{eq:affine-LPV-representation}.
\end{itemize}
Also, we discuss how such an approach can be realized for
flow control problems modelled by semi-discrete Navier-Stokes equations (NSE).
For this, we rely on
\begin{itemize}
  \item the coordinates provided by a proper orthogonal decomposition 
    (POD) of the velocity states; see~\cite{KunV02} for an introduction,
  \item the quadratic structure of the nonlinearity in the incompressible NSE,
  \item implicit treatment of the incompressibility con\-straint, and,
    importantly,
  \item low-rank solves for and low-rank representations of the solutions to the
    high-dimensional Riccati and Lyapunov equations.
\end{itemize}
We note that with this line of arguments, the feedback design by truncated
SDRE approximations can be made feasible for, say, finite element
approximations of general nonlinear partial differential equations (PDEs).

Apart from the proposed algorithmic advances and numerical insights into the
feedback approximation, we expand here on the work of~\cite{AllKS23} insofar as
the parametrization step lifts fundamental structural assumptions on the problem
class.
A related approach, though with updates that require the solutions of nonlinear
matrix equations, can be found in~\cite{CebC84} based on the expansion of
nonlinear systems into Volterra series~\cite{Rug81}.
Furthermore, we note that, with the explicit low-complexity parametrization of
the nonlinearity in an otherwise linear problem formulation, the difficulties
of exponentially growing dimensions that come with tensor expansions for general
nonlinearities are mitigated; see \cite{KraGB23} for a recent discussion
regarding model order reduction.

The overall procedure is explained in detail as follows.
In \Cref{sec:ld-aff-LPV}, we explain how a low-complexity linear
pa\-ram\-e\-ter-varying approximation can be obtained by pa\-ram\-e\-tri\-zing
the state of an SDC system.
The formulas for expanding the SDRE solution and state the constituting
equations for the coefficients of the expansion are recalled in
\Cref{sec:SDRE-expansion}.
\Cref{sec:nse-pod} lays out how POD can be used to realize a
low-complexity affine-linear parameter-varying approximation of incompressible
Navier-Stokes equations and in \Cref{sec:highdimme} we briefly describe
the concepts for solving the high-dimensional matrix equations.
In \Cref{sec:num-exps}, we provide the results of numerical experiments
to show the applicability of the approach and to compare with plain static
feedback.
The paper is concluded in \Cref{sec:conclusions}.


\section{Low-complexity linear parameter-varying approximations}
\label{sec:ld-aff-LPV}

In this section, we consider now systems in SDC form~\cref{eq:gen-sdc-sys}.
If the system state $v(t)$ is encoded into time-varying parameters
$\rho(t) = \mu(v(t)) \in \R^{m}$, with $m\leq n$, and
$v(t)= \nu(\rho(t))$, where $\mu$ and $\nu$ are the corresponding encoding and
decoding maps, then the
SDC representation~\cref{eq:gen-sdc-sys} can be formulated as a linear
parameter-varying (LPV) system via
\begin{equation} \label{eq:nonlsys-lpv}
  \dot{v}(t) = \tA(\rho(t)) v(t) + Bu(t), \quad y(t) = Cx(t),
\end{equation}
where $\tA(\rho) := A(\nu(\rho))$.
Such an embedding of a nonlinear system into the class of LPV systems is
typically called \emph{quasi LPV} system; see, e.g.,~\cite{KoeT20}.
Here we will focus on affine-linear LPV representations, where $\tA$ depends
affine-linearly on $\rho$ so that~\cref{eq:nonlsys-lpv} can be realized as
\begin{align*}
  \dot{v}(t)  = \left( \tA_0 + \sum_{k=1}^m \rho_{k}(t) \tA_{k} \right)
    v(t) + Bu(t), \quad
  y(t)  = C v(t),
\end{align*}
where $\rho_{k}$ is the $k$-th component of $\rho$ and where $\tA_{0}$ and
$\tA_k \in \R^{n\times n}$ are constant, $k=1,\ldots,m$.

If the state $v$ is not exactly parametrized but only approximated with less
degrees of freedom in $\hrho(t) = \hmu(v(t)) \in \R^{r}$ such that $r \ll m$,
with an inexact reconstruction
\begin{equation} \label{eq:approx-parametrization}
  v(t) \approx \tv(t) = \hnu (\hrho(t)) = \hnu (\hmu(v (t))),
\end{equation}
then an LPV approximation of~\cref{eq:gen-affine-system} is given by
\begin{equation} \label{eq:nonlsys-lpv-approx}
  \dot{\hv}(t) = \hA(\hrho(t)) \hv(t) + Bu(t), \quad \hy(t) = C \hv(t),
\end{equation}
with the approximated system matrix $\hA(\hrho) := A(\hnu(\hrho))$ and
the new system state $\hv(t)\in \R^{n}$.
Note that the state $\hv(t)$ is of \emph{full} dimension $n$, as our reduction
efforts will target the structure of the model rather than the dimension of the
states.
Nonetheless, the ideas and techniques of approximate low-dimensional
parametrizations of states and estimates on approximation errors of standard
model order reduction (MOR) schemes readily apply.


\section{Series expansions of state-dependent Riccati equations}%
\label{sec:SDRE-expansion}

The theory of first-order approximations for a single parameter dependency
by~\cite{BeeTB00} has been extended in~\cite{AllKS23} to the multivariate case.
We briefly recall the relevant formulas. 

To prepare the argument, we assume that $v$ is pa\-ram\-e\-tri\-zed through
$\rho$ and consider the dependency of the SDRE solution $P$ on the current
value of $\rho$, i.e., $P(\cdot) = P(\rho(\cdot))$.
Then, the multivariate Taylor expansion of $P$ about $\rho_{0} = 0$ up to order
$K$ reads 
\begin{equation} \label{eq:taylor-expansion-P}
  P(\rho) \approx P(0) + \sum_{1\leq |\alpha| \leq K} 
    \rho^{(\alpha)}P_{\alpha},
\end{equation}
where $\alpha =(\alpha_{1}, \ldots, \alpha_{m}) \in \N^{m}$ is a multiindex with
$|\alpha|:=\sum_{i=1}^{m} \alpha_{i}$, where
$\rho^{(\alpha)}:=\rho_1^{\alpha_1}\rho_2^{\alpha_2}\dotsm\rho_m^{\alpha_m}$,
and where, importantly, $P_\alpha$ are constant matrices given by
\begin{equation*}
  P_{\alpha} := \tfrac{1}{\alpha_{1}! \alpha_{2}! \cdots \alpha_{m}!}
    \tfrac{\partial^{|\alpha|}}
    {\partial_{\rho_{1}}^{\alpha_{1}} \partial_{\rho_{2}}^{\alpha_{2}}
    \cdots \partial_{\rho_{m}}^{\alpha_{m}}} P(0).
\end{equation*}
In particular, the expansion up to order one (i.e., the associated first-order
approximation) can be written as
\begin{equation} \label{eq:P-first-order}
  P(\rho) \approx P(0) + \sum_{|\alpha| = 1}  \rho^{(\alpha)}P_{\alpha} =: P_0 +
  \sum_{k=1}^m \rho_k L_k.
\end{equation}
Substituting $P$ in the SDRE~\cref{eq:sdre} by its series
expansion~\cref{eq:taylor-expansion-P} and considering the
affine-linear dependency of $A$ on $\rho$ yields
\begin{align*}
  & \bigl( \sum_{k=0}^{m} \rho_{k} A_{k} \bigr)^{\trans}
    \bigl( \sum_{|\alpha| \leq K} \rho^{(\alpha)} P_{\alpha} \bigr) +
    \bigl( \sum_{|\alpha| \leq K} \rho^{(\alpha)} P_{\alpha} \bigr)
    \bigl( \sum_{k=0}^{m} \rho_{k} A_{k} \bigr)\\
  & \quad{}-{}
    \bigl( \sum_{|\alpha| \leq K} \rho^{(\alpha)} P_{\alpha} \bigr)
    B B^{\trans}
    \bigl( \sum_{|\alpha| \leq K} \rho^{(\alpha)} P_{\alpha} \bigr)
    = -C^{\trans} C,
\end{align*}
where, for compactness of the expression, we introduce $\rho_0=1$ and use the
relevant conventions for $\alpha=0\in \mathbb N^{m}$.

By matching the coefficients for $\rho^{(\alpha)}$, we obtain equations for the
matrices of the first-order approximation~\cref{eq:P-first-order} as
\begin{equation} \label{eq:lpvsdre-base-riccati}
  A_{0}^{\trans} P_{0} + P_{0} A_{0} - P_{0} B B^{\trans} P_{0} = -C^{\trans} C,
\end{equation}
for $P_{0}$ and 
\begin{align} \nonumber
 & \left( A_{0} - B B^{\trans} P_{0} \right)^{\trans} L_{k} +
   L_{k} \left( A_{0} - B B^{\trans} P_{0} \right)\\
 \label{eq:lpvsdre-upd-lyapunovs}
 & = -\left(A_{k}^{\trans} P_{0} + P_{0} A_{k} \right)
\end{align}
for $L_{k}$, with $k = 1, \ldots, m$; see also~\cite{AllKS23}.


\section{Approximations of Navier-Stokes equations through linear affine
parametrizations}%
\label{sec:nse-pod}

As the standard model for incompressible flows we consider spatially discretized
Navier-Stokes equations (NSE).
After a shift of variables that eliminates constant nonzero Dirichlet boundary
conditions, the semi-discrete NSEs in the variables of the velocity
$v(t)\in \mathbb R^{n}$ and pressure $p\in \mathbb R^{n_p}$ with control input
$u$ can be written as
\begin{subequations} \label{eq:nse}
\begin{align}
  M\dot{v} &= \tN(v,v) + \tA v + J^{\trans} p + \tB u, \\
     0 & = J v. \label{eq:nse-conti-eq}
\end{align}
\end{subequations} 
At least for theoretical considerations, the incompressibility
constraint~\cref{eq:nse-conti-eq} can be resolved and the velocity $v$ can be
determined by the equivalent projected equations
\begin{equation} \label{eq:nse-projected}
    M \dot v = N(v,v) + Av  + Bu,
\end{equation}
where $N(v,v) = \Pi^{\trans} \tN(v,v)$, where $A$ and $B$ denote
$\Pi^{\trans} \tA$ and $\Pi^{\trans} \tB$, respectively, and where
\begin{equation} \label{eq:lerayproj}
  \Pi := I_{n} - M^{-1} J^{\trans}
    \left( J M^{-1} J^{\trans} \right)^{-1} J
\end{equation}
is the so-called \emph{discrete Leray projector}; see~\cite{morHeiSS08} for
properties of $\Pi$ and formulations in the coordinates of the
subspace spanned by $\Pi^{\trans}$ and see~\cite{BenH17b} where we have proven
that the SDRE feedback based on~\cref{eq:nse-projected} is equivalent to
that of~\cref{eq:nse}.

By the homogeneity in the boundary conditions, the nonlinearity $N(v,v)$ that
models the convection is linear in both arguments (see, e.g.,~\cite{BehBH17} for
explicit formulas of $N(\cdot,\cdot)$ in a spatial discretization) so that both
\begin{equation*}
  N_{1}(v) \colon w \mapsto N(v, w) \quad \text{and} \quad
  N_{2}(v) \colon w \mapsto N(w, v)
\end{equation*}
can be realized as state-dependent coefficient matrix and so that for any
blending parameter $\lambda \in \R$, an SDC representation is given as 
\begin{equation} \label{eq:nse-sdc}
  N(v,v) = \lambda N_{1}(v) v + (1-\lambda) N_{2}(v) v=: \Nl(v) v.
\end{equation}
Even more, if in an approximative parametrization as
in~\cref{eq:approx-parametrization}, the decoding is linear, then the induced
LPV approximation is \emph{affine-linear}; see~\cite[Rem. 2]{HeiBB22}.

In fact, let $V_{r} \in \R^{n\times r}$ be the matrix of $r$ POD modes
designed to best approximate the velocity in an $r$-di\-men\-sio\-nal subspace
of $\R^{n}$, then with
\begin{equation*}
  \hrho(t) := V_{r}^{\trans} v(t) \quad \text {and} \quad
  \tv(t) = V_{r} \hrho(t)
\end{equation*}
and the approximation property of the POD basis, we obtain
\begin{equation*}
  v(t) \approx \tv(t) = V_{r} \hrho(v(t)) =
    \sum_{i = 1}^{r} \hrho_{i}(v(t)) \hv_{i},
\end{equation*}
where $\hv_{i}$, for $i = 1, \dots, r$ are the columns of
$V_{r} \in \R^{n\times r}$.
By the linearity of $v\to \Nl(v)$, the SDC representation~\cref{eq:nse-sdc} is
readily approximated by
\begin{equation}\label{eq:pod-lpv-apprx}
  N(v) v \approx N(\tv) v = \left( \sum_{i = 1}^{r} \hrho_{i} \hN_{i} \right) v
\end{equation}
with $\hN_{i} := \Nl(\hv_{i})$, $i=1,\ldots,r$.

From the orthogonality of the POD basis it follows that
$\tv(t) = V_{r} V_{r}^{\trans} v(t) \to v(t)$ uniformly with $r \to n$,
$r\leq n$, which can be translated into convergence of the LPV approximations
\begin{equation*}
  \left( \sum_{i=1}^{\hr} \hrho_{i} \hN_{i} \right) 
  \to 
  \left( \sum_{i=1}^{r} \hrho_{i} \hN_{i} \right)
  \to 
  \left( \sum_{i=1}^{n} \hrho_{i} \hN_{i} \right) 
  = N(v),
\end{equation*}
for $\hr \to r \to n$ and $\hr \leq r \leq n$. 
Practically, in view of computing approximations to the SDRE solution and the
associated feedback law, this means that the series expansion in~\cref{eq:P-first-order}
can be augmented or reduced by simply adding or discarding parameters and
the corresponding factors $L_{k}$.


\section{Numerical handling of high-dimensional matrix equations}%
\label{sec:highdimme}

First, we note that for systems like the spatially discretized
NSEs~\cref{eq:nse-projected}, a mass matrix $M$ needs to be incorporated.
Since $M$ is typically symmetric positive definite and, thus, invertible, such
systems are readily transformed into the standard form
of~\cref{eq:gen-sdc-sys}.
In practice, however, it is beneficial to consider formulations of the Riccati
and Lyapunov equations~\cref{eq:lpvsdre-base-riccati}
and~\cref{eq:lpvsdre-upd-lyapunovs} without the explicit inversion:
\begin{equation} \label{eq:riccati_mass}
  A_{0}^{\trans} P_{0} M + M^{\trans} P_{0} A_{0} -
    M^{\trans} P_{0} B B^{\trans} P_{0} M = -C^{\trans} C,
\end{equation}
and
\begin{equation} \label{eq:lyap_mass}
  \Azcl^{\trans} L_{k} M + M^{\trans} L_{k} \Azcl =
    -\left( M^{\trans} P_{0} A_{k} +A_{k}^{\trans} P_{0} M \right),
\end{equation}
for $k = 1, \ldots, m$, where $\Azcl := A_{0} - B B^{\trans} P_{0} M$ is the
closed-loop system matrix corresponding to~\cref{eq:riccati_mass}; see,
for example,~\cite{BenH15}. 
Although, these formulations can cover also problems where the mass matrix is
not invertible as they appear in dif\-fer\-en\-tial-algebraic equations like the
Navier-Stokes equations in the original formulation~\cref{eq:nse}, for the
presentation it is convenient to refer to the projected
system~\cref{eq:nse-projected}.
In practice, the system matrices in~\cref{eq:nse-projected} and the involved
projection $\Pi$ from~\cref{eq:lerayproj} are realized only implicitly during
the application of iterative matrix equation solvers
for~\cref{eq:riccati_mass} and~\cref{eq:lyap_mass} like the low-rank ADI;
see~\cite{BenHW22}.

Another general problem occurring in the presence of high-dimensional systems
is the consumption of computational resources such as time and memory.
In particular, it is typically not possible to even store the large-scale dense
solutions $P_{0} \in \R^{n \times n}$ and $L_{k} \in \R^{n \times n}$,
$k = 1, \ldots, m$.
An established approach to handle this problem is the use of low-rank
factorizations.
The stabilizing solution of the Riccati equation~\cref{eq:riccati_mass} namely
$P_{0}$ is known to be positive semi-definite such that in the case of
small numbers of inputs and outputs, $p, q \ll n$, it can be well represented
by a low-rank Cholesky factorization, i.e.,
$P_{0} \approx Z_{0} Z_{0}^{\trans}$ with $Z_{0} \in \R^{n \times \ell_{0}}$ and
$\ell_{0} \ll n$.
One can find various numerical methods in the literature to efficiently compute
these low-rank factors of Riccati equations without ever forming the full
solution $P_{0}$; see~\cite{BenHW23} for an overview.
In our numerical experiments, we rely on the low-rank Newton-Kleinman-ADI
method~\cite{BenHSetal20}.

Considering the indefinite right-hand side of~\cref{eq:lyap_mass}, one needs to
assume that the $L_{k}$'s are generically indefinite.
Nonetheless, the low-rank factorization $P_{0} \approx Z_{0} Z_{0}^{\trans}$
yields the low-rank indefinite factorization
\begin{equation*}
  \begin{aligned}
    & M^{\trans} P_{0} A_{k} +A_{k}^{\trans} P_{0} M \\
    & \approx
    \begin{bmatrix} M^{\trans} Z_{0} & A_{k}^{\trans} Z_{0} \end{bmatrix}
    \begin{bmatrix} 0 & I_{\ell_{0}} \\ I_{\ell_{0}} & 0 \end{bmatrix}
    \begin{bmatrix} Z_{0}^{\trans} M \\ Z_{0}^{\trans} A_{k} \end{bmatrix}
    = C_{k} S_{k} C_{k}^{\trans}.
 \end{aligned}
\end{equation*}
Based on this factorization in the right-hand side, it is possible to similarly
approximate the solution to~\cref{eq:lyap_mass} as
$L_{k} \approx Z_{k} D_{k} Z_{k}^{\trans}$, for $k = 1, \ldots, m$, where
$D_{k}$ is a symmetric but possibly indefinite matrix.
Here, we are using the $LDL^{\trans}$-factorized
low-rank ADI method~\cite{LanMS14}.


\section{Numerical experiments}%
\label{sec:num-exps}

The code, raw data and results of the presented numerical experiments are
available at~\cite{supHeiW23}.
For the solution of Riccati and Lyapunov equations in MATLAB 9.9.0.1467703
(R2020b), we used the solver implementations from
MORLAB version 5.0~\cite{BenW19b} and
M-M.E.S.S. version 2.2~\cite{SaaKB22}.
For the simulation part, we resort to our Python interface \cite{swHei19} between \emph{Scipy} and the finite element toolbox FEniCS \cite{LogMW12}.

We consider the stabilization of the flow in the wake of a 2D cylinder
through two control inlets at the periphery of the cylinder.
Measurement outputs are defined as averaged velocities over a small
neighborhood of three sensor points in the wake; see~\cite{BehBH17} on technical
details, the detailed control setup, and on how the Dirichlet control is
relaxed as penalized Robin boundary conditions.

\begin{figure}
  \centering
  \includegraphics[width=\linewidth]{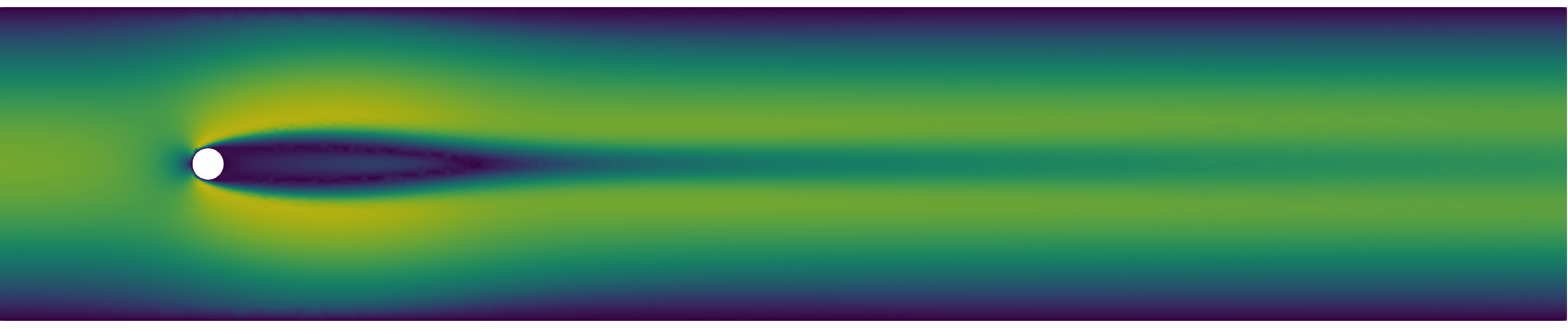}
  \includegraphics[width=\linewidth]{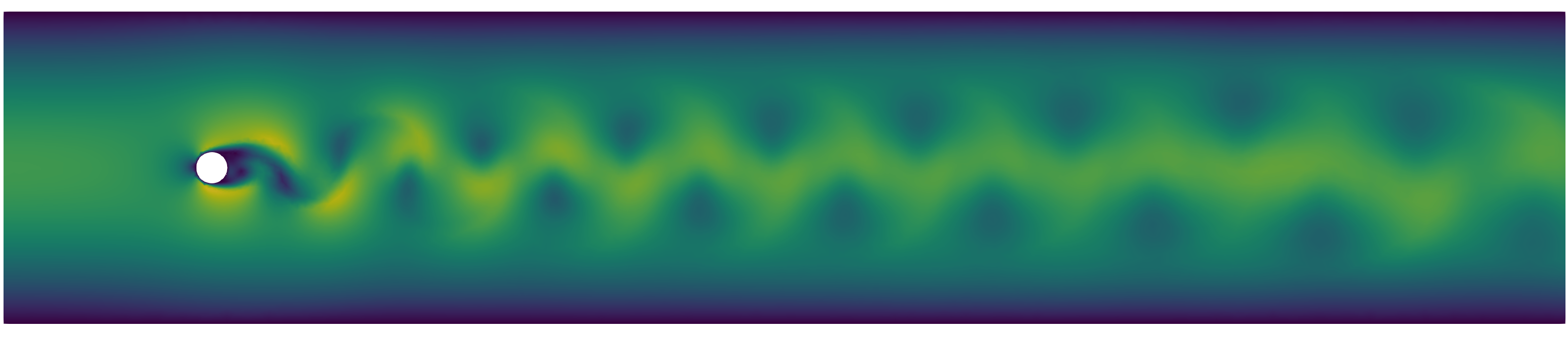}
  \caption{Snapshots of the flow in the unstable steady state and in the fully
    developed periodic vortex shedding regime.}
  \label{fig:Re60-flow-snapshots}
\end{figure}

As for the numerical setup, we consider here the Rey\-noldsnumber $60$ and start
from the associated non-zero steady state, which is to be stabilized;
see \Cref{fig:Re60-flow-snapshots} for the basic geometry of the example
and snapshots of the steady state and periodic regime that develops if no
stabilization is employed.
For the spatial discretization, we use quadratic-linear \emph{Taylor-Hood}
finite elements on a nonuniform mesh that leads to a system of size $57\,000$.
For the time integration we use an implicit-explicit Euler time stepping method
that in particular treats the linear part and the incompressibility constraint
implicitly, whereas the nonlinear part and the feedback is treated explicitly
in time.
Generally, we are concerned with a system of type~\cref{eq:nse} with
$\begin{bmatrix} v(t) & p(t) \end{bmatrix}^{\trans} \in \R^{57\,000}$ with the
input $u(t) \in \R^{2}$ and the output $y(t) \in \R^{6}$ extracted from the
velocity state by a linear output operator $C$ as $y(t) = C v(t)$.

The basic procedure of the simulations comprises the following steps:
\begin{itemize}
  \item[(0)] Compute the steady state for $u = 0$ to be used as reference for
    the stabilization, for the shift of the system that removes nonzero boundary
    conditions, and as the starting value for the 
    closed-loop simulations.
  \item[(1)] Perform open-loop simulations to collect data for the POD basis
    for the affine-linear LPV approximation~\cref{eq:pod-lpv-apprx}.
  \item[(2)] Compute the Riccati solution $P_{0}$ and the Lyapunov solutions
    $L_{k}$ via~\cref{eq:riccati_mass} and~\cref{eq:lyap_mass}.
  \item[(3)] Close the loop with the nonlinear feedback law
    \begin{equation}\label{eq:xsdre-fb}
      u(t) = -B^{\trans}
        \left( P_{0} + \sum_{k=1}^{r} \hrho_{k}(v(t)) L_{k} \right) M v(t).
    \end{equation}
\end{itemize}

In the presented numerical study, the relevant steps were realized as follows.
To acquire the data for the POD basis, we take $401$ snapshots of the velocity
equally distributed on the time interval $[0, 0.5]$ for the test signal 
\begin{equation} \label{eq:testu}
  u(t) = \begin{bmatrix} \sin(t) & 0 \end{bmatrix}^{\trans},
\end{equation}
and define $V_{r}$ as the matrix of the $r$ leading left singular vectors with
respect to the weighted inner-product induced by the mass matrix $M$
of the finite element discretization; see~\cite{BauBH18}.
Then the LPV approximation~\cref{eq:nse-sdc} is computed for the NSE
with $\lambda = 0.75$.
In the following, we denote the feedback definition by the truncated SDRE series approximation~\cref{eq:xsdre-fb} of parameter dimension $r$ by
$\xSDRE{r}$ and the classical linear-quadratic regulator feedback, which is
readily defined as $u(t) = -B^{\trans} P_{0} M v(t)$ by \LQR{}
(which is equivalent to \xSDRE{0}; cf. the feedback law~\cref{eq:xsdre-fb}).

To trigger the instabilities in the closed-loop simulations, we apply the test
signal~\cref{eq:testu} on a short time $[0, \tc]$ before we \emph{switch on}
the feedback at $t = \tc$.
In this way, the system will deviate from the linearization point.
For $\tc$ too large, the state may have left the region of attraction for which
the linear \LQR{} or the SDRE-based controller will stabilize the
nonlinear system.

\begin{figure}
  \centering
  \includegraphics[width=\linewidth]{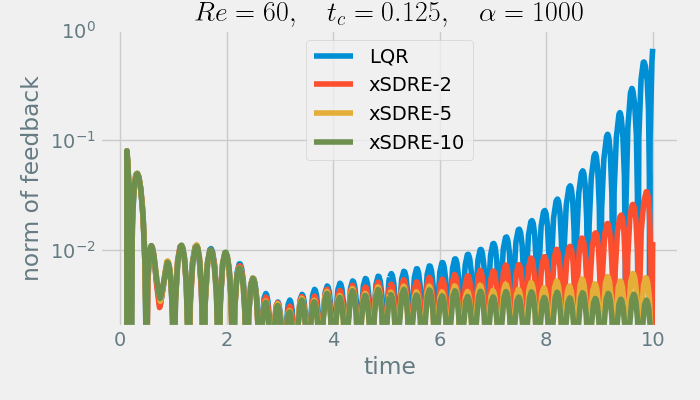}
  \includegraphics[width=\linewidth]{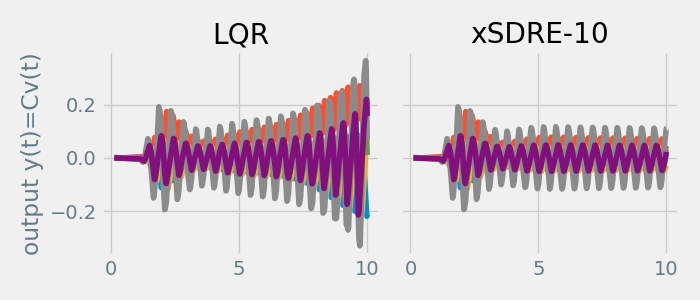}
  \caption{Norms of feedback for different $r$ and the corresponding outputs for
    the \LQR{} and the \xSDRE{10}-feedback for the case of $\alpha=10^3$ and
    $\tc=0.125$.
    A continuous performance improvement for \xSDRE{r} can be observed with
    stabilization of the system at \xSDRE{10}, while the classical
    \LQR{}-feedback fails.}
  \label{fig:re60-sut0125-fbs}
\end{figure}

In our experiments, we employed \emph{Tikhonov regularization} in the form of
an $\alpha \in \{ 1, 10^{3} \}$, which is included by replacing the original
input matrix $B$ by the scaled version $\breve{B} := \tfrac{1}{\sqrt{\alpha}} B$
in the definition of the SDRE~\cref{eq:sdre} as well as in the solved Riccati
and Lyapunov equations~\cref{eq:riccati_mass} and~\cref{eq:lyap_mass}.
Consequently, the corresponding feedback needs to be scaled as well, e.g.,
in the \LQR{} case as
\begin{equation*}
  u(t) = -\frac{1}{\alpha} B^{\trans} \breve{P}_{0} M v(t),
\end{equation*}
where $\breve{P}_{0}$ solves the Riccati equation with $\breve{B}$. 
Also, we used bisection of the time domain to identify a $\tc$ close to a
critical value that marks the performance region of the controllers. 
For both cases of $\alpha$, we found that the \xSDRE{r} approach enlarges the
domain of attraction of the \LQR-controller.
We illustrate this finding by plotting the outputs of the closed-loop systems as
measured for the \LQR-feedback and the \xSDRE{r}-feedback for various $r$.
Apart from the qualitative differences in the performance, e.g.,
stabilization is achieved or not, in these setups, a quantitative effect of the
reduction parameter $r$ as part of the nonlinear feedback definition becomes
evident.
In fact, for the case of $\alpha = 10^{3}$ and $\tc = 0.125$, the plain
\LQR-feedback did not stabilize the system, while additional
modes in the feedback design improved the performance steadily, until with
$r = 10$ stabilization was achieved;
see \Cref{fig:re60-sut0125-fbs} with the norms of the
feedback actions over time for this setup.

\begin{figure}
  \centering
  \includegraphics[width=\linewidth]{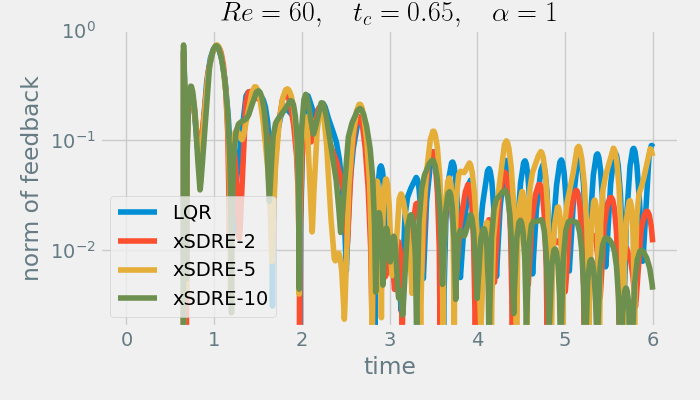}
  \includegraphics[width=\linewidth]{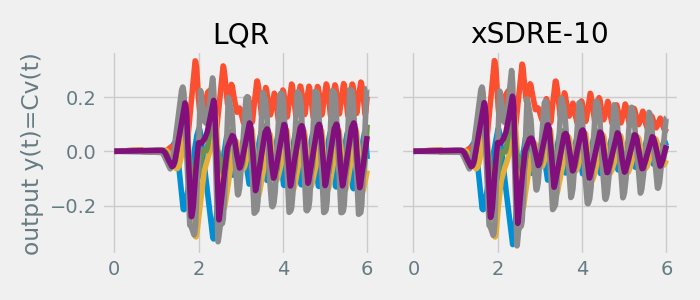}
  \caption{Norms of feedback for different $r$ and the corresponding outputs for
    the \LQR{} and the \xSDRE{10}-feedback for the case of $\alpha = 1$ and
    $\tc = 0.65$.
    The \xSDRE{10}-feedback is able to stabilize the system in contrast to the
    classical \LQR{} approach.}
  \label{fig:re60-sut0650-fbs}
\end{figure}

For smaller regularization parameter $\alpha = 1$ that enables larger control
actions, the overall region of performance is extended at the expense of a more
sensitive control regime.
Here, an illustrative $\tc$ could be identified at $t = 0.65$ with the
\LQR-feedback being not stabilizing in contrast to the
\xSDRE{10}-feedback; see \Cref{fig:re60-sut0650-fbs}.
However, with \xSDRE{2} was stabilizing while \xSDRE{5} failing to do so, a
clear trend for improvement in performance for larger $r$ could not be observed. 

Overall, we can state that the truncated SDRE approach, with an almost
negligible overhead in the simulation phase, is capable to improve the feedback
performance if compared to the classical LQR feedback.
In particular, for larger regularization parameters the difference in the
region of attraction was rather small but the smooth relation between the
additional complexity $r$ and the performance turned out to be a useful
hyperparameter in the feedback design. 
In the regime with less regularization, the \xSDRE{r}~approach proved to give
decisive advantages concerning the domain of performance at the extra cost of
finding a suitable $r$.


\section{Conclusion}%
\label{sec:conclusions}

In this work, we have presented a general framework that uses the embedding of
nonlinear systems in the class of LPV systems, POD for reduction of the
complexity of the parameter dimension, and the quadratic structure of the
convection term in the Navier-Stokes equations to make the nonlinear feedback
design through truncated expansions of the SDRE applicable. 
With state-of-the-art matrix equations solvers, computational feasibility of
this approach was achieved too. 
As illustrated by a numerical example, this generic nonlinear approach
provides a measurable improvement over the closely related classical LQR
approach with a small additional computational effort at runtime. 

Certainly, through the dependency on the POD basis for the parametrization, this
approach is problem-specific.
We would argue, however, that POD is a rather standard and generally applicable
way to define parametrizations so that our proposed method enjoys a similarly
general scope.

The potentials and needs for future work are manifold.
Firstly, it would be interesting to investigate expansions of second order.
Secondly, we have not paid particular attention to the definition of the
coordinates for the affine LPV approximation and simply resorted to POD. 
Since a small dimension $r$ is most important, in particular if one wants to
consider a second-order expansions of the SDRE, other, possibly nonlinear
parametrizations might be considered.
For the analysis of the approximation, suitable measures for the residual, e.g.,
in the SDRE approximation, need to be derived together with formulas for
feasible evaluations in the large-scale systems case.


\section*{Acknowledgments}%
\addcontentsline{toc}{section}{Acknowledgments}

Jan Heiland was supported by the German Research Foundation (DFG) Research
Training Group 2297 ``Mathematical Complexity Reduction (MathCoRe)'', Magdeburg.


\addcontentsline{toc}{section}{References}
\bibliographystyle{plainurl}
\bibliography{nse-lpv-sdre}

\end{document}